\documentclass[final]{siamltex}

\usepackage{amsfonts}

\begin{document}

\title{ON CONTROLLABILITY OF A NON-HOMOGENEOUS ELASTIC STRING WITH MEMORY}



\author{Sergei A. Avdonin \footnotemark[2] \footnotemark[4] \footnotemark[6]
\and Boris P. Belinskiy (corresponding author)\footnotemark[3] \footnotemark[5] \footnotemark[7]}
\maketitle

\renewcommand{\thefootnote}{\fnsymbol{footnote}}

\footnotetext[2]{University of Alaska Fairbanks, Fairbanks, AK 99775-6660}
\footnotetext[3]{University of Tennessee at Chattanooga, 615 McCallie Avenue, Chattanooga, TN 37403-2598}
\footnotetext[4]{s.avdonin@alaska.edu}
\footnotetext[5]{Boris-Belinskiy@utc.edu}
\footnotetext[6]{Research of this author was supported in part by the National Science
Foundation, grant ARC 0724860 }
\footnotetext[7]{Research of this author was supported in part by University of Tennessee at Chattanooga Faculty Research Grant}


\begin{abstract}
 We are motivated by the problem of control for a non-homogeneous elastic string with memory. We reduce the problem of controllability to a non-standard moment problem. The solution of the latter problem is based on an auxiliary Riesz basis property result for a family of functions quadratically close to the nonharmonic exponentials. This result requires the detailed analysis of an integro-differential equation and is of interest in itself for Function Theory. Controllability of the string implies observability of a dual system.
\end{abstract}


\def \b{\bigskip}
\def \beq{\begin{equation}}
\def \eeq{\end{equation}}
\def \la{\label}
\def \G{\Gamma}
\def \g{\gamma}
\def \de{\bigtriangleup}
\def \gr{\bigtriangledown}
\def \r{\rho}
\def \i{\infty}
\def \l{\lambda}
\def \sln{\sqrt{\lambda_n}}
\def \t{\tau}
\def \s{\sigma}
\def \p{\varphi}
\def \H{{\mathcal H}}
\def \L{{\mathcal L}}
\def \om{\omega}
\def \hom{\hat\omega}
\def \Om{\Omega}
\def \mN{\mathbb N}
\def \ent{e^{\pm i\om_n t}}
\def \ens{e^{\pm i\om_n s}}
\def \enta{e^{\pm i\om_n \tau}}
\def \en{e^{\pm i \om_n t+\nu t}}
\def \enm{e^{\pm i \om_n t-\mu t}}
\def \d{\delta}
\def \al{\alpha}
\def \be{\beta}
\def \th{\theta}
\def \k{\kappa}
\def \bu{\bullet}
\def \bom{\frac{b_n^\pm}{\omega_n}}
\def \ep{\epsilon}
\def \ET{\tilde \mathcal E}


\pagestyle{myheadings}
\thispagestyle{plain}
\markboth{Sergei A. Avdonin and Boris P. Belinskiy}{CONTROLLABILITY OF AN ELASTIC STRING WITH MEMORY}

\section{Introduction}
\la{sec:1} The conditions of controllability for the broad class of the linear oscillating structures have been under consideration since the classical papers by H.O. Fattorini and D.L. Russell (see \cite{Rus1}, \cite{DR}, the survey \cite{Rus3}, and book \cite{AI} for the history of the subject and extended list of references). The method of \cite{Rus1}, as well as many other subsequent papers, is based on the properties of exponential families (usually in the space $L^2(0,T)$), the most important of which for Control Theory are minimality and the Riesz basis property. Recent investigations into new classes of distributed systems such as hybrid systems and damped systems, as well as problems of simultaneous control, have raised a number of new and difficult problems in the theory of exponential families (see~\cite{AB}--\cite{AM01}, \cite{journal5},\cite{B},\cite{HZ}).

We say that a mechanical system, e.g. a string, is controllable if, for any initial data by suitable manipulation of the exterior forces, the system goes to the given regime. According to the classical scheme outlines above, the solution of the controllability problem is based on an auxiliary basis property result. The last result is of independent interest and represents the central part of the paper.

In this paper, we study controllability of an oscillating string the material of which has the {\sl memory}. Such a string is described by the equation of the form
\beq \la{PDE0}
 y_{tt}(x,t)=y_{xx}(x,t)+\int_0^t\,N(t-\t)y_{xx}(x,\t)d\t,\;(x,t)\in (0,l)\times(0,T),
\eeq
and we suppose that a control force is applied to one end of the string. Here and below, $N(\cdot)$ is the memory kernel.

In the short review of the literature below, we discuss the results on both heat and wave equation with memory since, formally speaking, controllability for the heat equation may be studied in a way similar to one for the wave equation. Also, the results on control for heat equations with memory motivated our study. We consider only linear models. Also we slightly change the authors' notations to make them uniform.

To the best of our knowledge, G.  Leugering has been the first author to study controllability for viscoelastic systems (see~\cite{Leug1,Leug2}). The model studied in \cite{Leug1} leads to the initial boundary value problem for the heat equation with a memory term. Since the case of the constant coefficients is considered, Laplace transform is used. Monotonicity of the memory kernel is essential for the proof of the exact controllability of the system in finite time. The existence of the time-optimal control is shown. An abstract differential equation of the second order with a memory term is considered in \cite{Leug2} and the generator of the appropriate semi--group is studied.

It has been known for a long time that the infinite speed of propagation exhibited by the heat equation is not quite physical. The corresponding modification of the Fourier law results in the heat equation with memory (see \cite{Catt}, \cite{GP}).

Controllability for the heat equation with memory is studied in \cite{BI} in the multidimensional case
\beq \la{heq1}
y_t(x,t)-\g\de y(x,t)-\int_0^t\,N(t-s)\de y(x,s)\,ds=f(x,t)\chi_{\om}(x)\,,
\eeq
$$(x,t)\in \Om\times(0,T),\;y(x,0)=y_0(x),\;x\in \Om;\;y(x,t)=0,\;(x,t)\in \partial \Om\times (0,T)\,.$$
Here the symbol $\chi_{\om}$ denotes the characteristic function of the open set $\om\subset \Om.$ The controllability is discussed under the strong assumptions on the memory kernel, which are alien to the mild assumptions we make in this paper. Specifically, $N(\cdot)$ is supposed to be {\sl completely monotone}, i.e. $(-1)^jN^{(j)}(t)\ge 0\;\forall t>0,\,j=0,1,...$. Under additional assumptions on its Laplace transform, approximate controllability is proved for $T>0$ if $\g >0.$ The one-dimensional equation on $x\in (0,l)$  is studied under the similar assumptions for the memory kernel. The exact controllability is proved for $T\ge l/\sqrt{N(0)}\,$ if $\g=0.$

In \cite{YZ}, a multi-dimensional problem with more general memory kernel is considered
\beq \la{heq2}
y_t(x,t)-\gr\cdot \int_0^t\,N(t-s,x)\gr y(x,s)\,ds=f(x,t)\chi_{\om}(x)\,,
\eeq
subject to the similar initial and boundary conditions. The memory kernel at $t=s$ is strictly positive (more precisely, $m_0\le N(0,x)\le m_1,\,\forall x\in \Om.$) Also some geometric conditions on the set $\om$ are assumed. Then there exists  $T_0>0$ such that the system is exactly controllable for all $T\ge T_0$ by means of control $f\in L^2[(0,T)\times \om].$ The method is based on Carleman estimates.

The paper \cite{YZ-a} represents a generalization of \cite{YZ} for anisotropic and nonhomogeneous medium, i.e. the equation
\beq \la{heq2-a}
y_t(x,t)-\sum_{i,j}\,\bigg(a^{ij}(x)\int_0^t\,N(t-s,x)y_{x_i}(x,s)\,ds\bigg)_{x_j}=f(x,t)\chi_{\om}(x)
\eeq
is considered. The similar results on exact controllability and observability are proved under the geometric conditions, which are more complex than in \cite{YZ}. On the other hand, if the geometric conditions are not satisfied, a negative controllability and observability result is presented.

If we consider (\ref{heq1}) with $\g=0,$ then differentiating (\ref{heq1}) with respect to time formally yields the wave equation with memory. Yet, generally speaking, controllability of the heat equation does not imply controllability of the wave equation.

The model
\beq \la{FaPo}
y_{tt}=N_0\de y+\int_0^t\,N(t-s)\de y(x,s)\,ds-a\,y_t
\eeq
is considered in \cite{FaPo}. Here $N_0>0,\,N\le 0,\,{\rm and}\,a\ge 0.$ The exponential decay of the solution is proved assuming the memory kernel $N(t)\not\equiv 0$ and decays exponentially.

The reachability problem for the wave equation with memory
\beq \la{Lor}
u_{tt}-u_{xx}+\b\int_0^t\,e^{-\eta(t-s)}u_{xx}(x,s)\,ds=0
\eeq
is studied in \cite{Lor}. We emphasize the special type of the memory kernel in (\ref{Lor}). The goal of the current paper is to prove controllability of a similar system with an arbitrary memory kernel (except  smoothness requirement).

In the next group of papers, the memory of the material is not necessarily taken  into consideration. Yet, these papers contain some technical ideas which are important to us. The wave equation with the time dependent tension is studied in \cite{AB}, \cite{AB-Bas}, \cite{ABI}, and \cite{ABP}. The proof of Riesz basis property of a family of the time-dependent functions in \cite{AB} and \cite{ABI} is based on the assumption that the tension varies ``slowly enough" with time. This restriction is removed in \cite{ABP}. To the best of our knowledge, the papers \cite{AB}, \cite{ABI}, \cite{ABP} represent the first attempt to apply the method of moments to equations with time dependent coefficients. The new difficulty in this case is the absence of an explicit representation for a family of functions arising in the moment problem. This fact, which  substantially complicates the analysis of controllability, is common to the heat equations with memory.

The papers \cite{P1,P3} are especially important to us.

The one--dimensional version of (\ref{heq1}) with $\g=0$ and without any monotonicity assumptions on $N$ is considered there:
\beq \la{Pa-1}
y_t(x,t)=ay(x,t)+\int_0^t\,N(t-\t)y_{xx}(x,\t)d\t\,.
\eeq
Controllability problem is reduced to a moment problem with respect to a family of functions. Techniques developed by L.~Pandolfi in \cite{P1,P3} allow proving that the aforementioned family forms a Riesz basis in a proper $L^2$ space, and that allows solving the exact controllability problem. In \cite{AP1,AP2} the problem of simultaneous temperature and flux controllability for heat equations with memory was studied.

Similarly to \cite{P1,P3,AP1,AP2}, we reduce the study of controllability to investigation of the Riesz basis property of the auxiliary functions, which we call ``quasi-exponentials''. However, there is a serious difference between the models discussed in these papers and our model. (a) In these papers, the transformation
\beq \la{mapsto}
y\mapsto e^{\th t}y\;\;{\rm with}\;\;\th=-N'(0)/N(0)
\eeq
leads to an equation of the same form as (1.7) but with $N'(0)=0.$ The last condition appears to be very helpful for the technical purposes. For the wave equation with memory, the similar simplifying condition has the form $N(0)=0$ which can not be achieved by a change of variables preserving the structure of equation. This fact sufficiently complicates the analytic problems we solve in the current paper. That is where the difference between our model and the models considered in  \cite{P1,P3,AP1,AP2} is essential. (b) We note also that the different techniques we use here results in the less restrictive requirements on the kernel $N(t)$ (see Assumption 2 below).

The present paper is the first work where the boundary controllability of a general 1d-wave equation with memory is studied. Sharp controllability time is established under mild regularity conditions for the (variable) coefficients of the equation.

The inclusion of the memory term in the model makes it necessary to analyze a non-standard moment problem associated with the solutions to a family of the Cauchy problems for an integro--differential equation. We use diverse methods of  Functional Analysis and Asymptotic Analysis to prove that these solutions (we call them quasi-exponential functions) possess the Riesz basis property in a proper $L^2$ space. This property allows us solving the moment problem and is of interest in itself for Function Theory.

The paper is organized as follows. In Section 2 we discuss the statement of the problem and the auxiliary Sturm--Liouville problem. In Section 3 we derive a series representation for the solution and the moment problem with respect to the quasi-exponential functions. Section 4 is the central in the paper. Here, we prove basis properties  results for these quasi-exponentials. In Section 5 we solve the moment problem and prove the main theorem about controllability, i.e. give the conditions of exact controllability of the string with memory. We also prove the observability of the dual system.

\section{Statement of the initial boundary value problem and controllability problem. The auxiliary Sturm--Liouville problem}
\la{sec:2} \setcounter{equation}{0} For any $T>0,$ we consider the following initial boundary value problem (IBVP)
 for the hyperbolic type PDE with the convolution term
\beq \la{PDE}
\;\;\r(x) y_{tt}(x,t)=(Ay)(x,t)+\int_0^t\,N(t-\t)(Ay)(x,\t)d\t,\;(x,t)\in (0,l)\times(0,T),
\eeq
\beq \la{bcond}
y_x(0,t)=f(t),\;y(l,t)=0\,,
\eeq
\beq \la{incond}
y(x,0)=y_t(x,0)=0\,.
\eeq
Here $A$ is the differential expression
\beq \la{A}
(A\p)(x) = (\al(x)\p'(x))'+\be(x)\p(x)\,.
\eeq
Above, $\r(x)$ is the density of the string; $l$ its length; $\al(x)$ its modulus of elasticity;
$N(\cdot)$ its memory kernel; the function $f(t)$ is control, i.e. the exterior force acting on the left end of the string.

The IBVP (\ref{PDE})--(\ref{incond}) describes the small (linear) transverse oscillations of the string with memory. Since the similar IBVP without memory term has a unique solution if the functions $\r,\,a,\,b,\,{\rm and}\,f$ are smooth enough (\cite{Ladyzh}), we may hope that the integral term in (\ref{PDE}) will not change this fact if the kernel $N$ is smooth enough.

This paper solves the exact controllability problem for the system (\ref{PDE})--(\ref{bcond}). More precisely, for any $T>0,$  we describe the reachability set
$$\bigg\{\big(y^f(\cdot,T),\,y^f_t(\cdot,T)\big):\;f\in L^2(0,T)\bigg\}$$
of the system.

We proceed further under the following

{\bf Assumptions.} {\sl

1. $\r$ and $\al$ are strictly positive functions from $C^2[0,l],$ $\be\in C[0,l]\,.$

2. $N''\in L^2(0,T)\,.$
}

{\bf Remark 1.} {\sl Assumption 1 is used to claim the well--known results on the auxiliary Sturm--Liouville problem (see Lemma 1 below). Using recent results on the Sturm--Liouville theory (see, e.g. \cite{VS}), we may weaken the regularity assumptions: $ \r'',\,a'',\,b\,\in L^1(0,l).$}

We further introduce the functional spaces,
\beq \la{sp}
\H\equiv L_{\r}^2(0,l),\;\;\H_1\equiv \{\p\in H^1(0,l),\,\p(l)=0,\}
\eeq
where $H^1(0,l)$ is the standard Sobolev space. Then there is a natural continuous embedding $\H_1 \subset \H,$ which leads to the natural embedding of $\H$ into the dual space $\H_{-1}:=\left( \H_1 \right)^\prime$ to $\H_1.$ We thus obtain the triple of spaces (so called {\em rigged space}):
$$\H_1 \subset \H \subset \H_{-1}\,.$$
We now may define the operator $A$ as the differential expression (\ref{A}) on the space $\H$ with the domain $ H^2(0,l) \cap \H_1\,.$

We consider the Sturm--Liouville problem associated with (\ref{A}), (\ref{bcond}):
\beq \la{St-L}
\begin{array}{c}
\displaystyle
(\al(x)\p'(x))'+\be(x)\p(x)+\l \r(x)\p(x)=0,\;0<x<l\,; \\
\displaystyle
\p'(0)=\p(l)=0\,.
\end{array}
\eeq

The following result is well--known (see \cite{Birkh}, Ch. X; \cite{WW}, Ch. VI; \cite{LevSar}; \cite{Mar}; \cite{N}; \cite{Fult}, and the references cited therein).

{\bf Proposition 1.} {\sl Let Assumption 1 be satisfied. Then the Sturm--Liouville problem (\ref{St-L}) has the following properties.

(a) Its spectrum $\{\l_n\}_{n=1}^\i$ is pure discrete, simple, and real with the only point of accumulation at $+\i.$

(b) The asymptotic representation holds,
\beq \la{la-as}
\om_n:=\sqrt{\l_n}=\sqrt{\l_n^0}+o(1)\;{\rm as}\;n\to \i\,.
\eeq
Here
\beq \la{T_0}
\sqrt{\l_n^0}=\frac{\pi}{L}\,\bigg(n-\frac{1}{2}\bigg),\;n\ge 1\,;\;\;
L=\int_0^l\,\sqrt{\frac{\r(x)}{\al(x)}}\,dx\,,
\eeq
($L$ is the optical length of the string, i.e. the time of propagation of the waves along the string). The eigenfrequencies $\om_n$ are separated,
\beq \la{om-sep}
\inf_{n\ne k}\,|\om_n-\om_k|>0\,.
\eeq

(c) The corresponding eigenfunctions  $\{\p_n\}_{n=1}^\i$ form an orthogonal basis in $\H$ (which we assume to be orthonormal),
\beq \la{orthog-1}
\int_0^l\,\r(x)\p_n(x)\p_m(x)\,dx=\d_n^m,\;\;(A\p_n,\p_m)_{L^2(0,l)}=\l_n\,\d_n^m\,.
\eeq

(d) The estimates hold
\beq \la{ph-n-as}
|\p_n(0)| \asymp 1.
\eeq
}
This relation means that
$$0 < \inf_{n \in \mN}  |\p_n(0)| \leq \sup_{n \in \mN}  |\p_n(0)| < \infty \,.$$

Proposition 1b implies

{\bf Proposition 2} \cite[Sec. II.4]{AI}, \cite{Seip}. {\sl Put $T_0=2L$. The following statements are valid.

(a) If $\l_n \ne 0$ for all $n,$ i.e. $\om_n \neq 0$ for all $n,$ then the family $\mathcal E=\{e^{\pm i\om_n \,t}\}_{n \in \mathbb{N}}$ forms a Riesz basis in $L^2(0,T_0).$

(b)  If $\l_n=0$ for some $n,$ then the statement (a) remains valid for the family
$\mathcal E=\{e^{\pm i\om_n\,t}\}_{n:\,\om_n\ne 0}\cup\{1\}\cup\{t\}\,.$

(c) For $T > T_0,$ there exists an infinite family of exponentials $\mathcal E_+^T$ such that the family $\mathcal E \bigcup \mathcal E_+^T$ forms a Riesz basis in $L^2(0,T)\,.$

(d) For $T < T_0,$ there exists an infinite subfamily  $\mathcal E_-^T \subset \mathcal E$ such that the family $\mathcal E \setminus \mathcal E_-^T$ forms a Riesz basis in $L^2(0,T)\,.$
}

All the statements of  Proposition 2 remain valid if we replace the family $ \mathcal E$ by $\{e^{\pm i\om_n \,t + \nu t}\}_{n \in \mathbb{N}}\,$
with any $\nu \in \mathbb C,$ since the operator of multiplication to $e^{\nu t}$ is bounded and boundedly invertible in $L^2(0,T)$ and therefore preserves the Riesz basis property.

{\bf Remark 2.} {\sl According to Proposition 1 (a), only a finite number of $\l_n$ may be negative. In this case, the corresponding $\om_n$ are pure imaginary, $\om_n=i\nu_n,$ and for definiteness, we assume $\nu_n>0\,.$}

The following result specifies the regularity of the solution to the IBVP (\ref{PDE})--(\ref{incond}).

{\bf Theorem 1.} {\sl Let Assumptions 1 and 2 be satisfied and $f\in L^2(0,T).$ Then $y\in C([0,T];\,\H_1)$ and $y_t\in C([0,T];\,\H)\,.$
}

Below we use the notation $y^f$ for the solution to the IBVP subject to the (controlling) force $f.$

The following results states the  controllability property of the IBVP (\ref{PDE})--(\ref{incond}).

{\bf Theorem 2.} {\sl Let Assumptions 1 and 2 be satisfied. Then

(a) For $T\ge T_0,$ the system (\ref{PDE})--(\ref{incond}) is exactly controllable, i.e.
\beq \la{Th-2}
\bigg\{\big(y^f(\cdot,T),\,y^f_t(\cdot,T)\big):\;f\in L^2(0,T)\bigg\}=\H_1\times \H\,.
\eeq

(b) For $T < T_0,$ the system (\ref{PDE})--(\ref{incond}) is not approximately controllable, i.e.
\beq \la{Th-2b}
{\rm cl} \ \bigg\{\big(y^f(\cdot,T),\,y^f_t(\cdot,T)\big):\;f\in L^2(0,T)\bigg\} \neq \H_1\times \H\,.
\eeq
Moreover,
\beq \la{Th-2c}
{\rm codim} \ {\rm cl} \ \bigg\{\big(y^f(\cdot,T),\,y^f_t(\cdot,T)\big):\;f\in L^2(0,T)\bigg\} \ \ {\rm in} \ \  \H_1\times \H\,= \infty\,.
\eeq
}

{\bf Remark 3.} {\sl The results of Theorems 1 and 2 can be easily extended to other boundary conditions and other types of boundary control and source-type controls like $g(x)f(t)$ with the given $g.$}

When this paper was ready to submission, we have been informed that the result close, though less general, to Theorem 2(a) (controllability of the system described by equation (\ref{PDE0}) for $T>2l$) was obtained in the paper by P. Loreti, L. Pandolfi and D. Sforza, ``Boundary controllability and observability of a viscoelastic string'' (submitted). 

\section{Moment problem}
\la{:3} \setcounter{equation}{0} As we show below, the moment problem originates from the terminal conditions
\beq \la{cont-y}
y(x,T)=v_0(x),\;y_t(x,T)=v_1(x),\;\;x\in [0,l],
\eeq
where
\beq \la{vv}
v_0 \in \H_1\,, \ \ v_1 \in \H\,.
\eeq
First, we derive an appropriate representation for $y(x,t).$ For that, we consider the inner product $(Ay,\p_n)_{L^2(0,l)}$ and integrate by parts twice to find
\beq \la{aux-ident-2}
\int_0^l\,((\al y')'+\be y)\p_n dx=\k_n f(t)-\l_n\,\int_0^l\,\r y\p_n\,dx\,.
\eeq
Here we use the  relations (\ref{bcond}) and (\ref{St-L}) and introduce the notation, $\k_n=-\al(0)\p_n(0).$ According to (\ref{ph-n-as}), $|\k_n|\asymp 1\,.$

We further expand the solution $y(x,t)$ to the IBVP (\ref{PDE})--(\ref{incond}) into the series with respect to the basis of the eigenfunctions $\p_n$ (see Proposition 1)
\beq \la{Four-1}
y(x,t)=\sum_{n\ge 1}\,a_n(t)\p_n(x)\,.
\eeq
We then multiply both sides of the basic Eq-n (\ref{PDE}) by an arbitrary eigenfunction $\p_n,$ integrate over $(0,l),$ and use the identity (\ref{aux-ident-2}). We finally find the integro--differential equation for the coefficients $a_n$(here and below we revert to the notation $\om_n=\sln\,$),
\beq \la{ode}
\ddot a_n(t)+\om_n^2\,a_n(t)+\om_n^2\int_0^t\,N(t-\t)a_n(\t)d\t=\k_n\,g(t)\,,
\eeq
subject to the homogeneous initial conditions (see (\ref{incond}))
\beq \la{incond-ode}
a_n(0)=0,\;\;\dot a_n(0)=0\,.
\eeq
Here we introduce the new ``control" $g(t),$
\beq \la{gf-rel}
g(t)=f(t)+\int_0^t N(t-\t)f(\t)\,d\t\,.
\eeq
We notice that the map $f \mapsto g$ is linear bounded and boundedly invertible in $L^2(0,T).$

The following estimate and equality  are standard (see, e.g. \cite[Sec. III.1]{AI})
\beq \la{est-1}
||y(\cdot,t)||_{\H_1}^2 \asymp \sum_{n\ge 1}\, n^{2} \, |a_n (t)|^2,\;\;||y_t(\cdot,t)||_{\H}^2 = \sum_{n\ge 1}\,|\dot a_n(t)|^2\,.
\eeq
Our next goal is to obtain a representation for $a_n(t).$ For that, we introduce two families $\{c_n(t)\},\,\{s_n(t)\},\,n\ge 1$ of solutions to the homogeneous integro--differential equation corresponding to (\ref{ode}), satisfying different initial conditions,
\beq \la{ode-c-s}
\begin{array}{c}
\displaystyle\ddot c_n(t)+\om_n^2\,c_n(t)+\om_n^2\int_0^t\,N(t-\t)c_n(\t)d\t=0,\;\;c_n(0)=1,\;\dot c_n(0)=0\,, \\
\displaystyle\quad\ddot s_n(t)+\om_n^2\,s_n(t)+\om_n^2\int_0^t\,N(t-\t)s_n(\t)d\t=0,\;\;s_n(0)=0,\;\dot s_n(0)=\om_n\,.
\end{array}
n\ge 1
\eeq
(In the absence of the memory, $N(t)\equiv 0,$ it follows $c_n(t)=\cos \om_n t,\;s_n(t)=\sin \om_n t.$) Then the solution of the Cauchy problem (\ref{ode-c-s}) has the form
\beq \la{a-n=}
a_n(t)=\k_n\,\int_0^t\,g(\t)\frac{s_n(t-\t)}{\om_n}\,d\t\,.
\eeq
From (\ref{ode-c-s}) it is easy to derive that $\dot s_n(t)/\om_n = c_n(t)$. Therefore,
\beq \la{a-n-dot}
\dot a_n(t)=\k_n\,\int_0^t\,g(\t)c_n(t-\t)\,d\t\,.
\eeq
Hence, we have found a representation (\ref{Four-1}), (\ref{a-n=}), (\ref{a-n-dot}) for the solution $y(x,t)$ and its time derivative $y_t(x,t).$ Then the  Eq-ns (\ref{cont-y}) can be written in the form
\beq \la{mom-I}
\sum_{n\ge 1}\,a_n(T)\p_n(x) = v_0(x),\;\; \sum_{n\ge 1}\,\dot a_n(T)\p_n(x) = v_1(x)\,.
\eeq
We further use the representation (\ref{a-n=}), (\ref{a-n-dot}) and the orthogonality condition (\ref{orthog-1}) to conclude
\beq
\la{asc}
\begin{array}{l}
\displaystyle\k_n\,\int_0^T\,g(\t)\,{s_n(T-\t)}\,d\t={\om_n} v_{0,n}, \\
\displaystyle\k_n\,\int_0^T\,g(\t)\,c_n(T-\t)\,d\t=v_{1,n},
\end{array}
\qquad
n\ge 1
\eeq
where we let
\beq \la{mom-aux}
v_{0,n}=(v_0,\p_n)_{\H},\;v_{1,n}=(v_1,\p_n)_{\H}\,.
\eeq
Inclusions (\ref{vv}) imply that $ \{\om_n v_{0,n}\}\,,\, \{v_{1,n}\} \in \ell^2.$

The equations (\ref{asc}) represent the desirable {\sl moment problem}: Find $g\in L^2(0,T)$ given $\{v_{0,n}\}\,,\,\{v_{1,n}\}\,.$

Theorem I.2.1 of \cite{AI} implies that the solvability of the moment problem (\ref{asc}) is the direct consequence of the basis property of the system $\{e_n^\pm(t)\}.$ First, we discuss the basis property (Section 4) and then solve the moment problem (Section 5).

Instead of the family $\{c_n\}\cup\{s_n\}$ we may consider the family $\ET=\{e_n^\pm(t)\}_1^\i,$
$$e_n^\pm(t)\equiv c_n(t)\pm is_n(t)$$
that satisfies the equations
\beq \la{e-n-pm}
\ddot e_n^\pm(t)+\om_n^2\,e_n^\pm(t)+\om_n^2\int_0^t\,N(t-\t)e_n^\pm(\t)d\t=0\,.
\eeq
The initial conditions for the functions $e_n^\pm(t)$ originate from the corresponding conditions for $c_n(t)$ and $s_n(t),$
\beq \la{incond-e-n-pm}
e_n^\pm(0)=1,\;\dot e_n^\pm(0)=\pm i\om_n\,.
\eeq
If for some $n,$ $\om_n=0,$ there is only one solution to the Cauchy problem (\ref{e-n-pm})--(\ref{incond-e-n-pm}), $e_n=1.$ In order not to ``lose" one solution, we proceed as follows. We note that the Eq-n (\ref{asc}) of the moment problem contains the ratio, $s_n (\cdot)/\om_n.$ The expression ${s_n(t)}/{\om_n}$ should be understood in the limiting sense if $\om_n=0.$ We show that for any finite $t,$
\beq \la{lim-om}
\lim_{\om_n\to 0}\,\frac{s_n(t)}{\om_n}=t\,.
\eeq
Indeed, the initial value problem (\ref{ode-c-s}) for $s_n(t)$ may be easily transformed into the following integral equation
\beq \la{lim-om-1}
s_n(t)=\sin\om_nt-\om_n\,\int_0^t\,\sin \om_n(t-\t)\int_0^\t\,N(\t-z)s_n(z)\,dz\,d\t\,.
\eeq
The limiting representation (\ref{lim-om}) may be now easily derived from (\ref{lim-om-1}). Therefore, if $\om_n=0$ for some $n,$ then we put the corresponding $e_n^\pm(t)$ in the family $\left\{e_n^\pm(t)\right\}$ to be $e_n^+=1$ and $e_n^-=t$ (see Proposition 2b).

\section{Basis property of $\ET=\{e_n^\pm\}$}
\la{:4} \setcounter{equation}{0}

In this section we develop new technical ideas that generalized the methods used in our previous papers \cite{AB,ABI} and in papers by L. Pandolfi \cite{P1,P3} (see also our joint paper \cite{ABP}).

We will use two definitions.

{\bf Definition 1}. {\sl A sequence $\{e_n\}$ is said to be $\om-$independent in a Hilbert space $H$ if  the conditions
$$\{\zeta_n\}\in l^2 \quad {\rm and}\quad \sum \zeta_n e_n=0$$
imply that $\zeta_n=0 $ for every $ n.$ The convergence of the series is understood in the norm of $ H.$
}

{\bf Definition 2}. {\sl A family $\{e_n\}$ is said to be an $\mathcal L$-basis in a Hilbert space $H$ if it forms a Riesz basis in the closure of its linear span.
}

The following theorem represents the main analytic result of this paper.

{\bf Theorem 3.} {\sl Let Assumptions 1 and 2 be satisfied. Then

(a) the family $\ET=\{e_n^\pm(t)\}$ forms a Riesz basis in $L^2(0,T_0)\,;$

(b) for $T > T_0,$ there exists an infinite sequence $\{{\tilde \om}_n\}$ of real numbers and the corresponding family ${\ET}_+^T$ of solutions to the Cauchy problems (\ref{e-n-pm})--(\ref{incond-e-n-pm})  such that the family $\ET \bigcup \ET_+^T$ forms a Riesz basis in $L^2(0,T)\,;$

(c) for $T < T_0,$ there exists an infinite subfamily  $\ET_-^T \subset \ET$ such that the family $\ET \setminus \ET_-^T$ forms a Riesz basis in $L^2(0,T)\,.$
}

Our proof is based on the Bari theorem \cite{Bari}, \cite[p. 317]{GK}, \cite[Theorem~15 p.~38]{Young}.

{\bf Theorem} (Bari). {\sl Any $\om-$linearly independent sequence $\{e_n\}$ which is quadratically close to some Riesz basis $\{\ep_n\}$ in a Hilbert space $H$, is itself a Riesz basis.
}

From Proposition 2 it follows that the family $\{e^{\pm i\om_n\,t+N(0) t/2}\}$ forms a Riesz basis in $L^2(0,2T_0)$ and an $\L$-basis in $L^2(0,T)$ for $T \geq 2T_0.$ Therefore, two steps suffice to prove the basis property of the family $\ET=\{e_n^\pm(t)\}.$

1) We find the asymptotic representation of the family $\ET=\{e_n^\pm\}$ to prove that it is quadratically close to the family
$\{e^{\pm i\om_n\,t+N(0) t/2}\}$ in $L^2(0,T).$

2) We prove that the functions $e_n^\pm(t)$ are $\om-$linear independent in $L^2(0,T)$ for $T \geq T_0.$

\subsection{ Quadratic closeness} At the first step, we may assume $t\in [0,T]$ with an arbitrary chosen $T>0.$ We construct the leading term of the asymptotic representation for $e_n^\pm(t)$ as $n\to \i.$ That allows us to prove the quadratic closeness. Our proofs are mostly based on numerous and cumbersome estimates with the main technical tool to be the Gronwall inequality along with the direct asymptotic estimates.

The integro--differential Eq-n (\ref{e-n-pm}) subject to the initial conditions (\ref{incond-e-n-pm}) is equivalent to the integral equation
\beq \la{e-n-1}
e_n^\pm(t)=\ent-\om_n\,\int_0^t\,\sin \om_n(t-\t)\,\int_0^\t\,N(\t-s)e_n^\pm(s)\,ds\,d\t
\eeq
$$=\ent-\om_n\,\int_0^t\,\sin \om_n z\int_0^{t-z}\,N(t-z-s)e_n^\pm(s)\,ds\,dz$$
or after integration by parts
\beq \la{e-n-22}
e_n^\pm(t)=\ent-\int_0^t\,N(t-s)e_n^\pm(s)\,ds+N(0)\int_0^t\,\cos \om_n z\,e_n^\pm(t-z)dz
\eeq
$$+\int_0^t\,\cos \om_n z\int_0^{t-z}\,N'(t-z-s)e_n^\pm(s)\,ds\,dz\,.$$
We rewrite the last integral equation as follows,
\beq \la{e-n-3}
e_n^\pm(t)=\ent+\int_0^t\,N_n^*(t-s)e_n^\pm(s)\,ds
\eeq
where the following kernel $N_n^*(\cdot)$ is introduced:
\beq \la{N-n-*}
N_n^{*}(t)\equiv -N(t)+N(0)\cos \om_n t+\int_0^t\,\cos \om_n z\,N'(t-z)dz\,.
\eeq
Below, we will need to analyze the asymptotic representation to the solution of a slightly more general integral equation than (\ref{e-n-3})
\beq \la{e-n-n-1}
v_n^\pm(t,\mu)=\ent+\int_0^t\,e^{-\mu(t-s)}\,N_n^*(t-s)\,v_n^\pm(s,\mu)\,ds
\eeq
with some parameter $\mu.$ Obviously $e_n^\pm(t)=v_n^\pm(t,0).$

Obviously the kernel $e^{-\nu t}\,N_n^*(t)$ is bounded provided $N'\in L^2(0,T),$
\beq \la{N-est}
\;\;\;\;\max_{[0,T]}\,|e^{-\mu t}\,N_n^*(t)|\le \big((1+c_1)\max_{[0,T]}\,|N(t)|+c_1\,\max_{[0,T]}\,\int_0^t\,|N'(z)|dz\big)e^{|\mu| T}
\equiv c_2,
\eeq
$$c_1=\max_{n}\,|\cos \om_n t|,\;t\in [0,T]\,.$$
(We remind here that some first $\om_n$ may be complex). More than that, $N_n^*\in C[0,T]\,.$
In (\ref{N-est}) and below, $c$ with a lower index denotes a positive constant that is independent of both $n$ and $t.$

Eq-n (\ref{e-n-n-1}) and the Gronwall inequality now imply
\beq \la{e-n-bound}
|v_n^\pm(t)|\le c_3,\;\;t\in [0,T]\,.
\eeq
This result, along with Eq-n (\ref{e-n-n-1}), implies that $v_n^\pm(t)\in C^2(0,T);$ for the function $\ent,$ this is also seen from (\ref{e-n-pm}).

To estimate $v_n^\pm(t)$ more precisely, we rewrite (\ref{e-n-n-1}) as follows
\beq \la{e-n-4}
v_n^\pm(t)=\ent+(A\,v_n^\pm)(t)+(B_n\,v_n^\pm)(t)+(C_n\,v_n^\pm)(t)
\eeq
where we introduce three integral operators
$$(A\,v_n^\pm)(t)\equiv-\int_0^t\,e^{-\mu(t-s)}\,N(t-s)v_n^\pm(s)\,ds;$$
$$(B_n\,v_n^\pm)(t)\equiv N(0)\int_0^t\,e^{-\mu(t-s)}\,\cos\om_n (t-s)\,v_n^\pm(s)dz;$$
$$(C_n\,v_n^\pm)(t)\equiv \int_0^t\,\cos \om_n z\int_0^{t-z}\,e^{-\mu(t-s)}\,N'(t-z-s)v_n^\pm(s)ds\,dz\,.$$

We proceed with the estimate of $v_n^\pm(t)$ in two steps.

{\bf Step 1}. {\sl The formal asymptotic representation for large $n.$}

First, we find
$$(C_n\,v_n^\pm)(t)=\int_0^t\,\frac{d\sin \om_n z}{\om_n}\,\int_0^{t-z}\,e^{-\mu(t-s)}\,N'(t-z-s)v_n^\pm(s)ds$$
$$=\int_0^t\,\frac{\sin \om_n z}{\om_n}\,\Big[N'(0)e^{\mu z}\,v_n^\pm(t-z)\,dz+\int_0^{t-z}\,e^{-\mu(t-s)}\,N''(t-z-s)v_n^\pm(s)\Big]\,dz\,,$$
so that
\beq \la{C-n-3}
|(C_n\,v_n^\pm)(t)|\le \frac{c_4}{n}\;\;{\rm if}\;\;N''\in L^2(0,T)\,.
\eeq

Everywhere below, the symbol $O$ means an estimate that is uniform with respect to $t \in [0,T].$

Since we proceed formally at this moment and are interested in the leading term of $v_n^\pm(t)$ as $n\to \i$ only, we remove the term $(C_n\,v_n^\pm)(t)=O\big(\frac{1}{n}\big)$ from the analysis, hence considering the (approximate) integral equation
\beq \la{e-n-5}
v_n^\pm=\ent+(A\,v_n^\pm)(t)+(B_n\,v_n^\pm)(t)\,.
\eeq
The solution of (\ref{e-n-5}) may be constructed by the (convergent) iteration series
\beq \la{iter-1}
v_n^\pm(t)=(I+A+A^2+...+B_n+B_n^2+AB_n+B_nA+...)\ent\,.
\eeq
We have $|A^k\,\ent|\le c_5/n.$ Indeed,
\beq \la{A-n}
(A\,\ent )(t)=-\int_0^t\,e^{-\nu(t-s)}\,N(t-s)\frac{d \ens}{\pm i\om_n}\,.
\eeq
Integrating by parts yields the desired estimate provided $N'\in L^1(0,T)$ (Assumption 2 is even stronger, $N''\in L^2(0,T).$) The estimate for the higher powers of $A$ is then obvious. We further evaluate the term $(B_n\,\ent)(t)$
\beq \la{B-n}
\begin{array}{c}
\displaystyle
(B_n\,\ent)(t)=N(0)\int_0^t\,e^{-\mu(t-s)}\,\cos \om_n(t-s)\ens\,ds= \\
\displaystyle
\frac{N(0)}{2\mu}\,\ent\,\big(1-e^{-\mu t}\big) +O\bigg(\frac{1}{n}\bigg)\,.
\end{array}
\eeq
Repeating the previous calculations yields
\beq \la{B-n-2}
(B_n^2\,\ent)(t)=\frac{N(0)}{2\mu}\,\enm\,\big(e^{\mu t}-1-\mu t\big)\,+O\bigg(\frac{1}{n}\bigg)\,.
\eeq
According to (\ref{A-n}) and (\ref{B-n}), the term $(AB_n\,\ent)(t)$ is of the order $O\big(\frac{1}{n}\big)$ (we again assume $N'\in L^2(0,T).)$

We conclude that the only terms in the iteration series (\ref{iter-1}) which produce the leading terms for $v_n^\pm(t)$ are generated by the leading terms of (\ref{B-n}), (\ref{B-n-2}), etc.
\beq \la{e-n-iter-1}
v_n^\pm(t)=(I+B_n+B_n^2+...)\ent\,.
\eeq

The leading term of the asymptotic representation for $(B_n^k\,\ent)(t),\,k> 2$ may be found with the help of the consecutive applications of the operator $B_n$ and removing the terms of the order $O\big(\frac{1}{n}\big).$ We finally get
\beq \la{e-n-iter-11}
v_n^\pm(t)=\sum_{j\ge 0}\,(B_n^j\,\ent)(t)=\enm\sum_{M\ge 0}\,\sum_{j\ge M}\,\frac{(\mu t)^j}{j!}\,\bigg(\frac{N(0)}{2\mu}\bigg)^M
\eeq
$$=\enm\,\sum_{j\ge 0}\,\frac{(\mu t)^j}{j!}\,\frac{(N(0)/2\mu)^{j+1}-1}{(N(0)/2\mu)-1}+...$$
$$=\enm\,\frac{1}{(N(0)/2\mu)-1}\Bigg(\frac{N(0)}{2\mu}\,e^{N(0)t/2}-e^{\mu t/2}\Bigg)+...$$
We need two specific cases of this formula.

{\bf Case 1}. If $\mu\to 0,$ the integral equation (\ref{e-n-n-1}) becomes the equation (\ref{e-n-3}). Hence, the representation for $e_n^\pm(0)t$ has the form
\beq \la{ent-as}
e_n^\pm(t)=\lim_{\mu\to 0}\,v_n^\pm(t,\mu)\equiv \en+...
\eeq
where we introduce the notation which appears below on the regular basis and will be used in the exponents only
\beq \la{nu}
\nu\equiv N(0)/2.
\eeq
{\bf Case 2}. If $\mu\to N(0)/2=\nu,$ the integral equation (\ref{e-n-3}) has the form
\beq \la{e-n-n-n-1}
v_n^\pm(t,\nu)=\ent+\int_0^t\,e^{-\nu(t-s)}\,N_n^*(t-s)\,v_n^\pm(s,\nu)\,ds
\eeq
and the following representation holds
\beq \la{v-n-as}
v_n^\pm(t)= \lim_{\mu\to \nu}\,v_n^\pm(t,\mu)=\ent(1+N(0)t/2)+...
\eeq

The correction terms in (\ref{ent-as}) and (\ref{v-n-as}) seem to be of the order $O\big(\frac{1}{n}\big)$ though that has not been proved yet.

We now briefly consider a particular case of the memory kernel $N(t-s)=a\,e^{-\eta (t-s)}$ (see \cite{Lor}) where the functions $e_n^\pm(t)$ may be found explicitly (and hence asymptotically). Indeed, if we apply Laplace transform to the problem
$$\ddot e_n^\pm(t)+\om_n^2\,e_n^\pm(t)+\om_n^2\int_0^t\,ae^{-\eta (t-s)}e_n^\pm(\t)d\t=0,$$
$$e_n^\pm(0)=1,\;\;\dot e_n^\pm(0)=\pm i\om_n,$$
we find
$$e_n^\pm(t)=\frac{1}{2\pi i}\int_C\,e^{pt}\,\frac{(p\pm i\om_n)(p+\eta)}{p^3+\eta p^2+\om_n^2 p+\om_n^2(a+\eta)}\,dp$$
with the usual contour of integration $C$ in the inverse Laplace transform. An elementary perturbation theory allows finding the poles of the integrand
$$p_{1,2}=\pm i\om_n+\frac{a}{2}+...,\;\;p_3=-a-\eta+...$$
Thus, the asymptotic representation for $e_n^\pm(t)$ is the same as in (\ref{e-n-iter-11})
$$e_n^\pm(t)=\en+O\bigg(\frac{1}{\om_n}\bigg)\,,$$
where we use the fact that $N(0)=a.$

{\bf Step 2}. {\sl Justification of the asymptotic representation (\ref{ent-as}).}

{\bf Lemma 1}. {\sl Let Assumptions 1 and 2 be satisfied. Then the family $\ET=\{e_n^\pm(t)\}$ is quadratically close in $L^2(0,T)$ to the family $\{\en\}.$}

{\bf Proof}. We actually need to prove that the order of the correction term in (\ref{ent-as}) is $O\big(\frac{1}{n}\big)$ uniformly on $t\in [0,T].$ For that, we introduce the new function,
\beq \la{c-n-2}
E_n^\pm(t)\equiv e_n^\pm(t)-\en\,.
\eeq
First of all, $E_n^\pm(t)\in C^2[0,T].$ We  now derive an integral equation for $E_n^\pm(t).$ For that, we substitute $e_n^\pm(t)=\en+E_n^\pm(t)$ into Eq-n (\ref{e-n-3}). After some long but elementary calculations, we find
\beq \la{E-n-1}
E_n^\pm(t)=\int_0^t\,N_n^*(t-s)E_n^\pm(s)\,ds+F_n^\pm(t),
\eeq
where we let
\beq \la{F-n-1}
F_n^\pm(t)=\underline{-\int_0^t\,N(t-s)e^{\pm i\om_n s+\nu s}\,ds}
\eeq
$$+\frac{N(0)}{N(0)\pm 4i \om_n}\,\Big(e^{\pm i\om_n t+\nu t}-e^{\mp i\om_n t}\Big)$$
$$+\underline{\int_0^t\,\cos \om_n (t-\t)}\int_0^\t\,N'(\t-s)e^{\pm i\om_n s+\nu s}ds\,d\t\,.$$
Integrating by parts in the underlined integrals yields
\beq \la{F-n-1-a}
F_n^\pm(t)=-\frac{N(0)}{\pm i\om_n+N(0)/2}\en+\frac{N(t)}{\pm i\om_n+N(0)/2}
\eeq
$$-\int_0^t\,N'(t-s)\frac{e^{\pm i\om_n s+\nu s}}{\pm i\om_n+N(0)/2}\,ds
+\frac{N(0)}{N(0)\pm 4i \om_n}\,\Big(\en-e^{\mp i\om_n t}\Big)$$
$$+\int_0^t\,\frac{\sin \om_n (t-\t)}{\om_n}\Big[\int_0^\t\,N''(\t-s)e^{\pm i\om_n s+\nu s}ds +N'(0)e^{\pm i\om_n \t+N(0)\t/2}\Big]\,d\t$$
so that
\beq \la{F-n-est-1}
F_n^\pm(t)=\frac{\phi_n^\pm(t)}{n}\;\;{\rm where}\;|\phi_n^\pm(t)|\le c_5\;\;{\rm if}\;N''\in L^2(0,T)\,.
\eeq
Applying the Gronwall inequality to Eq-n (\ref{E-n-1}) and using the estimates for the kernel (\ref{N-est}) and for the free term (\ref{F-n-est-1}) we find
\beq \la{E-n-est}
\max_{t\in [0,T]}\,|E_n^\pm(t)|\le \frac{c_6}{n}\,.
\eeq
Hence, according to (\ref{c-n-2}), the system $\{e_n^\pm(t)\}$ is quadratically close to the basis $\{\en\}\,.\;\;\bu$

The proof of $\om-$linear independence of $\{e_n^\pm(t)\}$ below requires analyzing the properties of  $(F_n^\pm(t))'$ and  $(F_n^\pm(t))''.$ Representation (\ref{F-n-1-a}) for $F_n^\pm(t)$ implies the uniform with respect to $t\in [0,T]$ asymptotic representation for its derivative
\beq \la{F-n-1-p}
(F_n^\pm(t))'=-N(0)\en+\frac{N(0)}{4}\,\Big(\en+e^{\mp i\om_n t}\Big)
\eeq
$$+\int_0^t\,\cos \om_n (t-\t)\int_0^\t\,N''(\t-s)e^{\pm i\om_n s+\nu s}ds\,d\t$$
$$+\frac{N'(0)}{N(0)}\,e^{\pm i\om_n t}\big(e^{\nu t}-1\big)+\psi_{1,n}(t)\,, \;\;{\rm where}\;\psi_{1,n}(t)=O\bigg(\frac{1}{\om_n}\bigg)\,.$$

Here we explicitly show all terms of the order $O(1)$ and introduce the correction term, $\psi_{1,n}(t).$ Representation (4.18) shows that if $\om_n=0$ for some $n,$ the corresponding $F_n^\pm$ is independent of $\om_n.$ In this case, we let $\psi_{1,n}=O(1).$ The similar representation for the second derivative has the form
\beq \la{F-n-1-pp}
(F_n^\pm(t))''=-N(0)(\pm i \om_n+\frac{N(0)}{2})\en
\eeq
$$+\frac{N(0)}{N(0)\pm 4i \om_n}\,\Big((\pm i \om_n+N(0)/2)^2\en-(i\om_n)^2e^{\mp i\om_n t}\Big)$$
$$-\om_n\,N'(0)\,\underline{\int_0^t\,\sin \om_n (t-\t)\int_0^\t\,N''(\t-s)e^{\pm i\om_n s+\nu s}\,ds\,d\t}$$
$$-\om_n\,N'(0)\underline{\int_0^t\,\sin \om_n(t-\t)e^{\pm i\om_n \t+N(0)\t/2}\,d\t}\,.$$

Before proving $\om-$linear independence of $\{e_n^\pm(t)\}$  we need some auxiliary results on the properties of the functions $F_n^\pm(t)$ and $E_n^\pm(t)\,.$

{\bf Lemma 2}. {\sl Let Assumptions 1 and 2 be satisfied. Then
\beq \la{F-n-0,F-n-1}
\begin{array}l
(a)\;\;F_n^\pm(0)=0,\;\,F_n^\pm(t)\in C^2[0,T],\\
(b)\;\;|(F_n^\pm)'(t)|\le c_7,\;\,|(F_n^\pm)''(t)|\le c_8\,\hom_n,\\
(c)\;\;N_n^*(0)=0,\;\,E_n^\pm(0)=0,\;\,|E_n^\pm(t)|\le \frac{c_6}{n}\le c_9,\\
(d)\;\;|(E_n^\pm)'(t)|\le c_{10},\;\,|(E_n^\pm)''(t)|\le c_{11}\hom_n\,.
\end{array}.
\eeq
where we let $\hom_n=|\om_n|$ if $\om_n\ne 0;\;\hom_n=1$ if $\om_n=0\,.$
}

{\bf Proof}. The results in (a) immediately follow from the representation (\ref{F-n-1}) and (\ref{F-n-1-pp}). The first two results in (c) immediately follow from (\ref{N-n-*}) and (\ref{E-n-1})   correspondingly. The third estimate in (c) is the direct consequence of (\ref{E-n-est}). The first estimate in (b) is a consequence of (\ref{F-n-1-p}). The second estimate in (b) is a consequence of (\ref{F-n-1-pp}). If $\om_n=0$ for some $n,$ the corresponding $(F_n^\pm)(t)$ is independent of $\om_n$ (see (\ref{F-n-1})). To prove (d), we differentiate Eq-n (\ref{E-n-1}) and integrate by parts to find
\beq \la{E-n'-est}
(E_n^\pm)'(t)=\int_0^t\,(N_n^*)'(t-s)E_n^\pm(s)\,ds+N_n^*(0)E_n^\pm(t)+(F_n^\pm)'(t)
\eeq
$$=-N_n^*(0)E_n^\pm(t)+N_n^*(t)E_n^\pm(0)+\int_0^t\,N_n^*(t-s)(E_n^\pm)'(s)\,ds+(F_n^\pm)'(t)\,.$$
Since $N_n^*(0)=E_n^\pm(0)=0$ (see (\ref{F-n-0,F-n-1}) (c)) we find the final form of the equation for $(E_n^\pm)'(t)$ to be
\beq \la{E-n'-eq}
(E_n^\pm)'(t)=\int_0^t\,N_n^*(t-s)(E_n^\pm)'(s)\,ds+(F_n^\pm)'(t)\,.
\eeq
The unique solvability of this Volterra integral equation in $C[0,T]$ follows from the continuity of the free term $(F_n^\pm)'(t)\,,$ see Lemma 2 (a).] Applying the Gronwall inequality to Eq-n (\ref{E-n'-eq}) and an estimate of $|(F_n^\pm)'(t)|$ from (\ref{F-n-0,F-n-1}) (b) yield
\beq \la{E-n'-est-1}
|(E_n^\pm)'(t)|\le c_{10}\,.
\eeq
Differentiating (\ref{E-n'-eq}) yields
\beq \la{E-n''-eq}
(E_n^\pm)''(t)=\int_0^t\,(N_n^*)'(t-s)(E_n^\pm)'(s)\,ds+(F_n^\pm)''(t)
\eeq
$$=-N_n^*(t-s)(E_n^\pm)'(s)|_{s=0}^{s=t}+\int_0^t\,N_n^*(t-s)(E_n^\pm)''(s)\,ds+(F_n^\pm)''(t)\,.$$
Since $N_n^*(0)=0$ (see Lemma 2 (c)) and $(E_n^\pm)'(0)=(F_n^\pm)'(0)$ (which follows from (\ref{E-n'-eq})) we find
\beq \la{E-n''-est-1}
(E_n^\pm)''(t)=N_n^*(t)(F_n^\pm)'(0)+\int_0^t\,N_n^*(t-s)(E_n^\pm)''(s)\,ds+(F_n^\pm)''(t)\,.
\eeq
The unique solvability of this Volterra integral equation in $C[0,T]$ follows from the continuity of the free term $(F_n^\pm)''(t)\,,$ see Lemma 2 (a).]  Applying the Gronwall inequality to Eq-n (\ref{E-n''-est-1}) and using the estimate of $|(F_n^\pm)''(t)|$ (see  (\ref{F-n-0,F-n-1}) (b)), we find
\beq \la{E-n''-est-2}
|(E_n^\pm)''(t)|\le c_{11}\hom_n\,.\;\;\bu
\eeq

\subsection{ $\om-$ linear independence}
\verb+ +

We proceed in two steps, which are similar but not identical to ones in \cite{ABP,P3} (see also the discussion in Introduction).

The next result describes the asymptotic structure of the sequence $\{(E_n^\pm)'(t)\}.$ Lemma 2(d) only shows that it is bounded.

{\bf Lemma 3.} {\sl The sequence $\{(E_n^\pm)'(t)\}$ admits the following asymptotic representation
\beq \la{E-1}
\;\;\;\;\;\;\;\;(E_n^\pm)'(t)=D_1\,\en+D_2\,\en\cdot N(0)t/2+q_{n,1}(t)
\eeq
where the sequence $\{q_{n,1}(t)\}\in l^2$ uniformly on $[0,T]$ and the constants $D_1$ and $D_2$ are independent of $n.$\,}

{\bf Proof}. The free term of the integral equation (\ref{E-n'-eq}) is the linear combination of the following terms,
\beq \la{RHS-1}
\ent,\,\en,\,\psi_{1,n}(t),
\eeq
and the integral term
\beq \la{RHS-2}
I_n(t)\equiv\int_0^t\,\cos \om_n (t-\t)\int_0^\t\,N''(\t-s)e^{\pm i\om_n s+\nu s}ds\,d\t.\,
\eeq
Here the sequence $\{\psi_{1,n}\}\in l^2$ uniformly on $t\in [0,T].$ In the integral $I_n(t),$ we substitute $\t-s\mapsto s$ and change the order of integration. After that, the interior integral may be evaluated. We find
\beq \la{RHS-3}
\;\;\;\;\;I_n(t)=\int_0^t\,N''(s)e^{\mp i\om_ns-\nu s}\,\Bigg[\ent\,\frac{e^{\nu t}-e^{\nu s}}{N(0)}+
O\bigg(\frac{1}{\om_n}\bigg)\Bigg]\,ds
\eeq
where the last $O$ is uniform on $[0,T].$ Assumption $N''\in L^2(0,T)$ implies that $I_n(t)\in l^2$ uniformly on $(0,T).$ We conclude that the free term is the linear combination of three terms, a sequence from $l^2$ uniformly on $(0,T)$ (see (\ref{RHS-2}) and (\ref{RHS-3})) and two exponentials (see (\ref{RHS-1})). Hence the solution $(E_n^\pm)'(t)$ may be represented as a linear combination of the  corresponding solutions, which we denote as $u_{n,j}^\pm(t),\,j=1,2,3.$ We denote the proofs of the estimates for $u_{n,j}$ as {\bf Steps E1, E2,} and {\bf E3}.

{\bf Step E1}. First, we show that if the free term of the integral equation (\ref{E-n'-eq}) belongs to $l^2$ uniformly on $[0,T],$ then the solution has the same property. Indeed, the solution $u_{n,1}(t)$ of the integral equation
$$u_{n,1}^\pm(t)=\int_0^t\,N_n^*(t-s)u_{n,1}^\pm(s)\,ds+r_n(t),\;\{r_n(t)\}\in l^2,\,$$
may be estimated as follows
$$|u_{n,1}^\pm(t)|\le |r_n(t)|+c_1\,\int_0^t\,|u_{n,1}^\pm(s)|\,ds.$$
According to the Gronwall inequality
$$|u_{n,1}^\pm(t)|\le |r_n(t)|\,e^{c_1t}\le e^{c_1T}\,|r_n(t)|\;{\rm so\;that}\;\{u_{n,1}^\pm(t)\}\in l^2.$$

{\bf Step E2}. The integral equation with $\ent$ in the free term coincides with the equation (\ref{e-n-3}). Hence
$$u_{n,2}^\pm(t)=\en+O\Big(\frac{1}{n}\Big).\,$$

{\bf Step E3}. We use the results of Lemma 1 here. We consider the integral equation
\beq \la{int-E-1}
u_{n,3}^\pm(t)=\en+\int_0^t\,N_n^*(t-s)\,u_{n,3}^\pm(s)\,ds\,,
\eeq
rewrite it as
$$u_{n,3}^\pm(t)e^{-\nu t}=\ent\,\int_0^t\,N_n^*(t-s)\,e^{-\nu (t-s)}\,u_{n,3}^\pm(s)e^{-\nu s}\,ds\,,$$
and introduce the new unknown function, $v_n^\pm(t)\equiv u_{n,3}^\pm(t)e^{-\nu t},$ which satisfies the equation (\ref{e-n-n-n-1}). Its asymptotic representation for $v_n^\pm(t)$ is given by Lemma 1, see (\ref{v-n-as}). We conclude
\beq \la{u-n-3-as}
u_{n,3}^\pm(t)=\en\,(1+N(0)t/2)+...
\eeq
The justification of the last asymptotic representation follows one in Lemma 1. We introduce the function
\beq \la{int-E-3}
U_n(t)^\pm\equiv u_{n,3}^\pm(t)-\en\,\Big(1+N(0)t/2\Big)
\eeq
which satisfies the integral equation similar to (\ref{e-n-3}) or (\ref{E-n-1})
\beq \la{int-E-4}
U_n^\pm(t)=\int_0^t\,N_n^*(t-s)\,U_n^\pm(s)\,ds+G_n^\pm(t).
\eeq
For the free term $G_n^\pm(t),$ after some long but elementary calculations, we find
\beq \la{int-E-5}
G_n^\pm(t)=\en-\en\,\Big(1+N(0)t/2\Big)+
\eeq
$$N(0)\,\int_0^t\,\cos \om_n(t-s)\,e^{\pm i\om_n s+\nu s}\,\Big(1+N(0)s/2\Big)\,ds
+O\bigg(\frac{1}{n}\bigg)=O\bigg(\frac{1}{n}\bigg)\,.$$
Repeating Step 1 yields $\{U_n(t)\}\in l^2.$ This shows that representation (\ref{u-n-3-as}) actually holds with the error term from $l^2,$ and this completes the proof of representation (\ref{E-1}). $\bu$

The next result describes the asymptotic structure of the sequence $\{(E_n^\pm)''(t)\}.$ Lemma 2(d) only shows that it may not grow faster than $\hat\om_n.$

{\bf Lemma 4.} {\sl The sequence $\{(E_n^\pm)''(t)\}$ admits the following asymptotic representation
\beq \la{E-2}
\;\;\;\;\;\;\;\;(E_n^\pm)''(t)=\hat\om_n\Big(D_4\,\en+D_5\,e^{\mp i\om_n t+\nu t}+D_6\,\ent\,N(0)t/2+q_{n,1}(t)\Big)
\eeq
where the sequence $\{q_{n,2}(t)\}\in l^2$ uniformly on $[0,T]$ and the constants $D_{*}$ are independent of $n.$\,}

{\bf Proof}. The proof is similar to the proof of Lemma 3. We consider the integral equation (\ref{E-n''-est-1}) and analyze the free term in it, $N_n^*(t)(F_n^\pm)'(0)+(F_n^\pm)''(t).$ The direct calculation shows that $(F_n^\pm)'(0)=-N(0)/2$ and also
\beq \la{F-n-2-as}
(F_n^\pm)''(t)=\pm i\om_n\,N(0)\Big(-\frac{3}{4}\,\en-\frac{1}{4}\,e^{\mp i\om_n t}\pm
\eeq
$$i\,\frac{N'(0)}{N(0)}\,\int_0^t\,\sin \om_n (t-\t)\int_0^\t\,N''(\t-s)e^{\pm i\om_n s+\nu s}ds\,d\t\Big)\,.$$
As in the proof of Lemma 3, we conclude that the free term in (\ref{E-n''-est-1}) is the linear combination of the following terms,
\beq \la{RHS-4}
\en,\,e^{\mp i\om_nt},\,\psi_{2,n}(t),
\eeq
where the sequence $\{\psi_{2,n}\}\in l^2$ uniformly on $t\in [0,T].$ The remaining part of the proof is the same as in Lemma 3. $\bu$

{\bf Lemma 5}. {\sl Let Assumptions in Section 2 be satisfied and $T\ge T_0.$ Then the family $\ET=\{e_n^{\pm}(t)\}$ is $\om-$independent in $L^2(0,T)\,.$}

{\bf Proof}. The proof consists of several steps.

{\bf Step A} {\bf (Linear independence.)} We consider an arbitrary finite subfamily $\{e_n^{\pm}(t)\}_{n\le n_0}$ and demonstrate that the equality
\beq \la{om-ind-fin}
\sum_{\pm}\sum_{n=1}^{n_0}\,a_n^\pm e_n^\pm(t)=0
\eeq
is not possible for a nontrivial finite sequence $\{a_n\}_{n\le n_0}.$ Differentiating the representation (\ref{e-n-1}) for $e_n^\pm(t)$ implies (we note that $\l_n=\om_n^2$)
$$\big(e_n^\pm(t)\big)'=\pm i\om_n e^{\pm i\om_n t} +\l_n\,\int_0^t\,\cos \om_n(t-\t)\,\int_0^\t\,N(\t-s)e_n^\pm(s)ds\,d\t\,;$$
\beq \la{e-n''}
\big(e_n^\pm(t)\big)''=-\l_ne^{\pm i\om_n t}-\om_n^3\int_0^t\,\sin \om_n (t-\t)\,\int_0^\t\,N(\t-s)e_n^\pm(s)ds\,d\t
\eeq
$$+\l_n\int_0^t\,N(t-\t)e_n^\pm(\t)d\t=-\l_n\Big(e_n^\pm(t)-\int_0^t\,N(t-s)e_n^\pm(s)ds\Big)\,.$$
Using the last representation and differentiating the identity (\ref{om-ind-fin}), we find
\beq \la{e-n''-1}
-\sum_{\pm}\sum_{n=1}^{n_0}\,a_n^\pm\l_ne_n^\pm(t)+\big(N*\sum_{\pm}\sum_{n=1}^{n_0}\,a_n\l_ne_n^\pm\big)(t)=0\,.
\eeq
The operator $-I+N*$ is invertible, which implies
\beq \la{om-ind-2-a}
\sum_{\pm}\sum_{n=1}^{n_0}\,a_n^\pm\l_ne_n^\pm(t)=0\,.
\eeq
Transition from (\ref{om-ind-fin}) to (\ref{om-ind-2-a}) shows that similarly
\beq \la{om-ind-2-b}
\sum_{\pm}\sum_{n=1}^{n_0}\,a_n^\pm\l_n^ke_n^\pm(t)=0,\;k=0,1,...,n_0-1\,.
\eeq
The notation
$$\eta_n=\sum_{\pm}\,a_n^\pm\,e_n^\pm(t),\;1\le n\le M$$
allows to view the series of conditions (\ref{om-ind-2-b}) as a linear algebraic system with Vandermonde matrix $\{\l_n^k\},\;0\le k\le n_0-1,\;1\le n\le n_0\,.$ Hence, all $\eta_n=0\,.$ Linear independence of the functions $e_n^+$ and $e_n^-$ (for every index $n$) implies $a_n^\pm=0\,.$ Yet, the finiteness of the sum ($n_0<\i$) is required for  this conclusion.

{\bf Step B}. We prove that the equality
\beq \la{om-ind-2}
\sum_{\pm}\sum_{n\ge 1}\,a_n^\pm e_n^\pm(t)=0
\eeq
is not possible for a nontrivial sequence $\{a_n^\pm\}\in l^2,$ where convergence is understood in the $L^2(0,T)-$norm; $T\ge T_0\,.$

We introduce the function
\beq \la{E-n-11}
E(t)=-\sum_{\pm}\sum_{n\ge 1}\,a_n^\pm e^{\pm i\om_n t+\nu t}\,.
\eeq
Since $\{a_n^\pm\}\in l^2,$ $E \in L^2(0,T).$

{\bf Step B1}. We prove that $a_n=b_n/\om_n,$ where $\{b_n\}\in l^2\,.$ It is assumed here that none of $\om_n=0,\,$ and this assumption is not restrictive. For example, we may divide $b_n$ by $n$ instead of $\om_n.$

According to (\ref{c-n-2}) and (\ref{om-ind-2}),
\beq \la{e-cap-n-2}
E(t)=\sum_{\pm}\sum_{n\ge 1}\,a_n^\pm\,E_n^\pm(t)\,.
\eeq
Using $\{a_n^\pm\}\in l^2$ yields $\sum\,|a_n^\pm|/n<\i\,.$ Since $|E_n^\pm(t)|\le c_6/n$ (see (\ref{E-n-est})) we conclude by the Weierstrass $M-$test that $E(t)\in C[0,T].$  Also by Lemma 2(c),   $E_n^\pm(0)=0$ for every $n,$ so that $E(0)=0\,.$ We further compute the derivative of $E(t)$ given by (\ref{e-cap-n-2}) term-wise and then prove that the resulting series converges in $L^2(0,T).$ We find
\beq \la{e-cap-n-3}
E'(t) \sim \sum_{\pm}\sum_{n\ge 1}\,a_n^\pm\,(E_n^\pm)'(t)\,.
\eeq
Here and below the symbol $\sim$ denotes the formal differentiation of a series. We study this series with the help of representation (\ref{E-1}) (see Lemma 3). Substituting the term $q_{n,1}(t)$ into the series yields a uniformly convergent series since the sequence $\{a_n\}\in l^2.$ Substituting any of the remaining exponentials produces a series that converges in $L^2(0,T)$ since the set of exponentials forms an $\L$-basis in $L^2(0,T),\,T\ge T_0.$

Hence, $E(t)$ given by (\ref{E-n-11}) or (\ref{e-cap-n-2}) belongs to $H^1(0,T),$ so that $$\sum_{\pm}\,\sum_{n\ge 1}\,a_n^\pm\,e^{\pm i\om_n t}\in H^1(0,T)\,.$$
Since the family $\{e^{\pm i\om_n t}\}$ forms an $\L$-basis in $L^2(0,T),\,T \ge T_0,$ the family $\{\ent/\om_n\}$ forms an ${\cal L}-$basis in $ H^1(0,T)$ \cite{Rus2}, \cite[Sec. II.5]{AI}. We conclude that $\{a_n^\pm\om_n\}\in l^2$ or $a_n^\pm=b_n^\pm/\om_n,$ where $\{b_n^\pm\}\in l^2\,.$

{\bf Step B2}. We prove that $a_n=c_n/\l_n,$ where $\{c_n\}\in l^2$ (we remind here that $\l_n=\om_n^2$). Again, differentiating the series (\ref{e-cap-n-2}) twice (formally) yields
\beq \la{E''-1}
E''(t)\sim \sum_{\pm}\sum_{n\ge 1}\,\bom(E_n^\pm)''(t)\,.
\eeq
We study this series with the help of representation (\ref{E-2}) (see Lemma 4), according to which the structure of the ratio $(E_n^\pm)''(t)/\om_n$ is the same as the structure of $(E_n^\pm)'(t).$ We conclude that the series (\ref{E-1}) converges in $L^2(0,T).$

We conclude that $E \in H^2(0,T).$ Representation
(\ref{E-n-11}) now shows that
\beq \la{E-in-H-2}
E(t)=-\sum_{\pm}\sum_{n\ge 1}\,a_n^\pm e^{\pm i\om_n t+\nu t}\in H^2(0,T)\,.
\eeq
Hence, $\{a_n^\pm\l_n\}\in l^2$ or $a_n^\pm=c_n/\l_n,$ where $\{c_n\}\in l^2\,.$ Since we may now differentiate twice the series (\ref{om-ind-2}) term-wise, we can repeat the proof from {\bf Step A} above. We conclude that all $a_n^\pm=0,$ and this completes the proof of $\om-$independence.$\;\;\bu$

Analyzing the proofs in Sections 4.1 and 4.2, we see that our results  can be formulated as follows.

{\bf Theorem 4.} {\sl Let the family $\left\{e^{ i\mu_n t}\right\}$ forms an $\L$-basis in $L^2(0,T),$
and functions $e_n(t)$ satisfy the equations
\beq \la{e-mu}
\ddot e_n(t)+\mu_n^2\,e_n(t)+\mu_n^2\int_0^t\,N(t-\t)e_n(\t)d\t=0
\eeq
with the initial conditions
\beq \la{emu}
e_n(0)=1,\;\dot e_n(0)= i\mu_n\,.
\eeq
Then the family  $\left\{e_n (t)\right\}$ is also an $\L$-basis in $L^2(0,T)\,.$
}

To prove this statement we notice that the $\L$-basis property of $\left\{e^{ i\mu_n t}\right\}$ implies that
\beq \la{mun}
\sup \Im |\mu_n| < \i\,, \  \ |\mu_n| + 1 \geq C\cdot n
\eeq
with some positive constant $C$ (see, e.g. \cite[Secs. II.1, II.4]{AI}; the addition of $1$ in the last inequality reflects the possibility of $\mu_n=0$ for one single $n.$) In turn, the proofs in Section 4.1 demonstrate that conditions (\ref{mun}) imply the quadratic closeness of the families $\left\{e_n (t)\right\}$ and $\left\{e^{ i\mu_n t + N(0) t/2}\right\}\,.$

Theorem 3 follows now from Theorem 4 and Proposition 2.

\section{Controllability. Observability}
\la{:5} \setcounter{equation}{0}
Now we apply the basis property results (Theorem 3) to the moment problem (\ref{asc}) which we may rewrite after simple change of notations in terms of the family $\{e_n^\pm(t)\},\,n\ge 1$ as follows
\beq \la{a-n-mom-1}
a_n^{\pm}=\k_n\,\int_0^T\,g(t)e_n^\pm(t)\,dt\,
\eeq
with $\{a_n^{\pm}\} \in \ell^2\,.$

Theorems 1 and 2 follow now from the results of \cite{AI} (Theorems I.2.1(a,e), III.3.10(a), and Lemma III.2.4). We can summarize the facts we use from the cited references as follows.

{\sl Since the family $\left\{e_n^\pm (t)\right\}$ forms an $\L$-basis in $L^2(0,T)$ for $T \geq T_0,$ then

(i) if  $T \geq T_0,$ the moment problem (\ref{a-n-mom-1}) is solvable for any $\{a_n^{\pm}\} \in \ell^2$;

(ii) for any $T>0,$ the function $t \mapsto \left\{ \int_0^tg(\t)e_n^\pm(\t)d\t \right\}$ is continuous from $[0,T]$ to $\ell^2.$

Since
for $T < T_0,$  the family $\ET \setminus \ET_-^T$ forms a Riesz basis in $L^2(0,T)\,,$ then  the map $g \mapsto \{a_n^{\pm}\}$ from
$L^2(0,T)$ to $\ell^2$ defined by (\ref{a-n-mom-1}) has the range with codimension equal to the number of elements in $\ET_-^T$. }

As usually, the controllability result can be presented as observability of the dual system. Using the standard techniques based on the integration by parts, one can check that the dual system to (\ref{PDE})--(\ref{incond}) can be written in the form (after the change of variables $T-t \to t$)
\beq \la{PDE-d}
\;\;\;\r(x)v_{tt}(x,t)=(Av)(x,t)+\int_0^t\,N(t-\t)(Av)(x,\t)d\t,\;(x,t)\in (0,l) \times (0,T),
\eeq
\beq \la{bcond-d}
v_x(0,t)=v(l,t)=0,
\eeq
\beq \la{incond-d}
v(x,0)=v_0(x),\;v_t(x,0)=v_1(x); \ v_0 \in \H, \; v_1 \in \H_{-1}\,.
\eeq
The equivalent form of Theorems 1 and 2 reads as follows.

{\bf Theorem 1a.} {\sl Let Assumptions 1 and 2 be satisfied. Then for any $T>0$ the estimate
\beq \la{emb-v-1}
||v(0,\cdot)||_{L^2(0,T)}^2\le C \,\big(||v_0||_{\H}^2+||v_1||_{\H_{-1}}^2\big)
\eeq
holds with a positive constant $C$ independent of $v_0,\,v_1.$
}

{\bf Theorem 2a.} {\sl Let Assumptions 1 and 2 be satisfied. Then for $T\ge T_0 $ the estimate
\beq \la{emb-v-2}
||v(0,\cdot)||_{L^2(0,T)}^2\ge C \,\big(||v_0||_{\H}^2+||v_1||_{\H_{-1}}^2\big)
\eeq
holds with a positive constant $C$ independent of $v_0,\,v_1.$
}

Equivalence of Theorems 1 and 1a and, correspondingly, Theorems 2 and 2a expresses the standard relations between a linear operator and its adjoint in a Hilbert space. Specifically for partial differential equations, relations between Theorem 1 and 1a are described by the {\it transposition method} \cite{LM}. The estimate (\ref{emb-v-2}) presents the {\it observability inequality} for the system (\ref{PDE-d})--(\ref{incond-d}) which is equivalent to exact controllability of (\ref{PDE})--(\ref{incond}). More details about applications of the duality principle to hyperbolic equations with memory can be found in \cite{Lor}.



\section{Acknowledgments}

The authors are grateful to Luciano Pandolfi for fruitful discussions.
\vskip 1cm


\end{document}